\documentclass{mymathart}
\usepackage{mymathmacros}
\newtheorem{thm}{Theorem}[section]

\newtheorem{prop}[thm]{Proposition}
\newtheorem{lem}[thm]{Lemma}
\newtheorem{cor}[thm]{Corollary}
\newtheorem{defn}[thm]{Definition}

\usepackage{hyperref}

\newcommand{\Ek}{E_{\kappa}}

\newcommand{\extaddress}[1]{\underline{#1}}

\newcommand{\Sequ}{\mathcal{S}}

\newcommand{\Sequb}{\overline{\Sequ}}

\renewcommand{\u}{{\tt u}}
\newcommand{\m}{{\tt m}}

\newcommand{\bdyit}[2]
             {{\rule{0pt}{0pt}_{\mbox{$\scriptstyle #2$}}^{\mbox{%
                   $\scriptstyle #1$}} }}
\renewcommand{\j}{{\tt j}}
\newcommand{\itj}{\bdyit{\j}{\j-1}}

\newcommand{\addr}{\extaddress{r}}

\newcommand{\ADDR}{\operatorname{Addr}}
\newcommand{\Addr}{\operatorname{Addr}}
                               
\newcommand{\addu}{\extaddress{\u}}
\newcommand{\ut}{\tilde{\u}}
\newcommand{\addut}{\extaddress{\ut}}

\newcommand{\adds}{\extaddress{s}}
\newcommand{\s}{\adds}

\newcommand{\addt}{\extaddress{t}}
\renewcommand{\t}{\addt}

\newcommand{\extaddr}{\operatorname{addr}}
\newcommand{\itin}{\operatorname{itin}}
\newcommand{\K}{\mathbb{K}}

\newcommand{\gr}{g_{\addr}}

\newcommand{\gs}{g_{\adds}}

\newcommand{\ts}{t_{\s}}

\renewcommand{\H}{\mathbb{H}}

\renewcommand{\r}{\extaddress{r}}

\title{On Nonlanding Dynamic Rays of Exponential Maps}
\author{Lasse Rempe}
\address{Dept.~of Mathematical Sciences, University of Liverpool, 
  Liverpool L69 7ZL,
United Kingdom}
\email{l.rempe@liverpool.ac.uk}
\date{\today}
\thanks{Supported in part
 by a postdoctoral fellowship of the 
 German Academic Exchange Service (DAAD)}


\begin{document}

 \begin{abstract}
  We consider the case of an exponential map
  $\Ek:z\mapsto\exp(z)+\kappa$ for which the singular value $\kappa$
  is accessible from the set of escaping points of $\Ek$. We show that
  there are dynamic rays of $\Ek$ which do not land. In particular,
  there is no analog of Douady's ``pinched disk model'' for exponential
  maps whose singular value belongs to the Julia set. 
 We also prove that
  the boundary of a Siegel disk $U$ for which the singular value is accessible
  both from the set of escaping points and from $U$ contains 
  uncountably many indecomposable continua. 
 \end{abstract}
 \maketitle

 \section{Introduction}

In polynomial dynamics, 
 \emph{dynamic rays} 
 provide an important 
 tool which permits the investigation of
 a function's dynamics in 
 combinatorial terms. In many
 important cases, the Julia set is locally 
 connected, and all dynamic
 rays land. In this situation, the 
 Julia set can be described 
 as a 
 ``pinched disk'' \cite{pincheddisk}, 
 that is, as the quotient of $\mathbb{S}^1$ by a
 natural equivalence relation. These
 ideas form the foundation for many
 spectacular advances recently made in
 the study of polynomial maps.

In this article, we consider the family of
 exponential maps
 $\Ek:z\mapsto \exp(z)+\kappa$, which has 
 enjoyed much attention
 over the past two
 decades as the simplest parameter space
 of transcendental entire functions. 
 For such maps, the 
 set of escaping points 
  \[ I(\Ek) := \{z\in\C:|\Ek^n(z)|\to\infty\} \]
 contains no open subsets, and thus belongs to the Julia set, which
 is not locally connected even in the simplest cases. 

 Nonetheless, it is known that the 
  path-connected components of
  $I(\Ek)$ are curves to $\infty$ 
  \cite{expescaping}, providing an
  analog for dynamic rays of polynomials. 
  It is therefore natural to ask whether 
  there is also 
  an analog of the notion of
  local connectivity, and a corresponding
  topological model for the dynamics of
  ``tame'' exponential maps defined
  only in terms of their combinatorics.

This is indeed the case when
 $\Ek$ has an attracting or 
 parabolic cycle: here the Julia 
 set is homeomorphic to a ``pinched 
 Cantor Bouquet''; that is, the quotient
 of a certain (universal) space by
 a suitable equivalence relation. Furthermore,
 the dynamics on $J(\Ek)$ is the quotient of
 a universal dynamical system on this Cantor Bouquet. 
 (See \cite[Theorem 5.7]{accessible} and
   \cite[Corollary 9.3]{topescapingnew}.)
 We will show that
 this is the only case in which this
 is possible. 

\begin{thm}[Exponential Maps with Simple Julia Sets] \label{thm:nopincheddisk}
 Let $\Ek$ be an exponential map.
  Then the following are equivalent:
\begin{enumerate}
 \item $\Ek$ has an attracting
    or parabolic orbit.
  \item Every dynamic ray of $\Ek$
    lands in $\Ch$ and every point
    of the Julia set is on a
    dynamic ray or the landing
    point of such a ray.
\end{enumerate}
\end{thm}

This shows that there is no
 obvious analog of the 
 ``pinched disk model'' for
 exponential maps whose singular
 value $\kappa$ belongs to the
 Julia set. This includes ``tame'' 
 examples such as the
 postsingularly finite (``Misiurewicz'')
 case. In 
 contrast, polynomial Misiurewicz maps
 always have locally connected Julia sets.

To prove Theorem 
 \ref{thm:nopincheddisk}, we will
 show the following
 result on the existence of
 nonlanding rays.

\begin{thm}[Existence of Nonlanding Rays] \label{thm:nonlanding}
 Suppose that $\Ek$ is an exponential 
 map whose singular value
  $\kappa$ is on a dynamic ray
  or is the landing point of such a ray.
 
  Then there exist uncountably many 
  dynamic rays $g$ whose
  accumulation set (on the Riemann 
  sphere) is an indecomposable continuum
  containing $g$.
\end{thm}

If the hypotheses of this theorem
 are satisfied, we say that the
 singular value is \emph{accessible}
 (from the escaping set).
 By \cite{expescaping,expper}, 
  all Misiurewicz parameters, as well
  as
  all parameters for 
  which the singular
  value
  escapes (\emph{escaping parameters}) 
  satisfy
  this 
  condition. We expect that
  accessibility holds for a much larger
  class of parameters.

The presence of indecomposable 
 continua in exponential dynamics 
 (albeit not as the accumulation
  set of a dynamic ray) was 
 first
 observed by Devaney 
 \cite{devaneyknaster} 
 when
 $\kappa\in (-1,\infty)$.  For 
 Misiurewicz
 parameters, the existence of dynamic
 rays whose accumulation sets are such 
 continua (as in Theorem
 \ref{thm:nonlanding}) was first 
 observed by
 Schleicher in 2000 (personal 
 communication; see also \cite{misindecomposable}).
 The same result was proved
 in the case of 
 $\kappa\in (-1,\infty)$ by Devaney and
 Jarque \cite{indecomposable}, using
 similar methods.
 The basic idea underlying both these
 results, as well as our proof,
 is fairly simple (see ``idea of the 
 proof'' below).
 However, we will gain control 
 of the accumulation sets of
 dynamic rays using combinatorial
 considerations rather than the 
 previously used 
 expansion methods 
 (which rely on the fact that 
 the singular orbit is discrete).
 This allows us to prove our
 result in a 
 considerably 
 more general 
 situation.

Theorem \ref{thm:nonlanding} has the 
 interesting continuum-theoretic
 corollary that there is no analog of the 
 Moore triod theorem \cite[Proposition 2.18]{pommerenke} for
 ``Knaster-like'' continua; i.e.\ 
 indecomposable continua
 containing dense 
 curves (compare \cite[Section 8]{pommerenkecarmonacontinua}
 for a discussion of such questions). 

\begin{cor}[Uncountable number of indecomposable continua] 
 \label{cor:uncountablymanydisjointcontinua}
  There is an uncountable set of pairwise disjoint indecomposable
  plane continua each of which contains a dense injective curve. 
\end{cor}

\begin{figure}
 \fbox{%
  \resizebox{.48\textwidth}{!}{%
    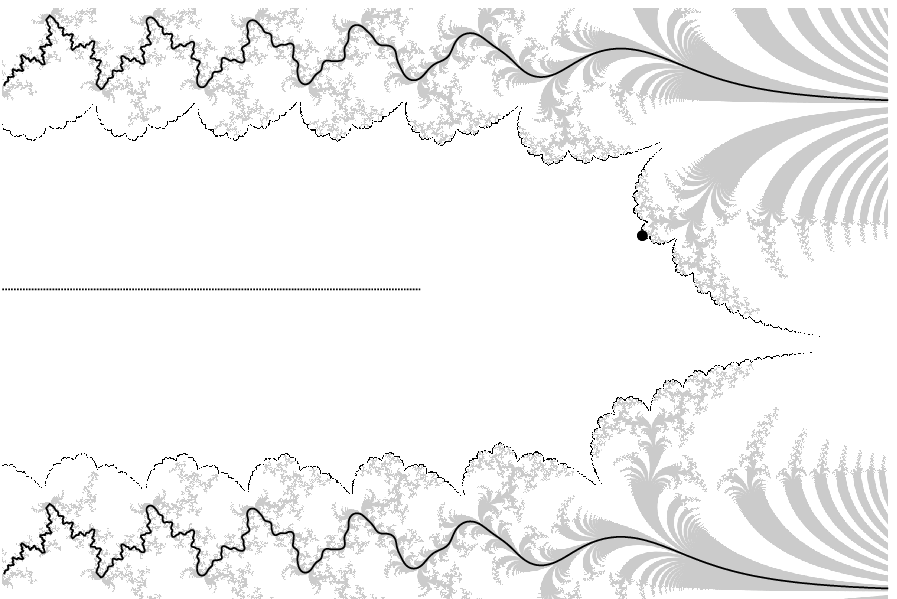}}%
   \hfill%
 \fbox{%
  \resizebox{.48\textwidth}{!}{%
    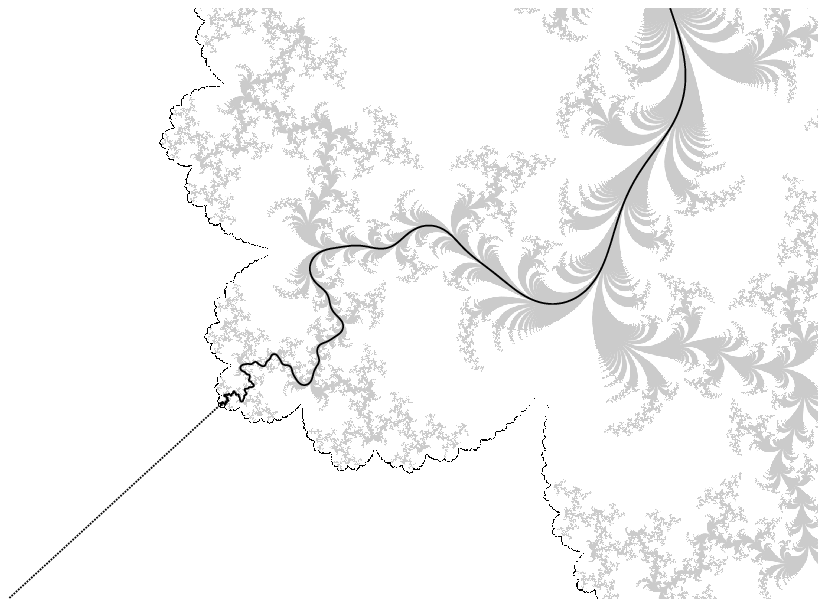}}%
   \hfill%
 \caption{An exponential map with a Siegel disk $U$ whose rotation number is
   the golden mean. The solid curves are dynamic rays 
   (contained in the escaping set $I(\Ek)$), while the dotted line is in
   the Siegel disk. The singular orbit (which is dense in $\partial U$) is also
   plotted.
   The picture on the right is a magnification about the singular value; the
    curves here are the images of the corresponding curves on the left. These
    pictures suggest that $\kappa$ is accessible both from $U$ and from
    $I(\Ek)$.\label{fig:siegel}}
\end{figure}

 Our results also have an interesting consequence for exponential
  maps with Siegel disks. Suppose that $\Ek$ has a Siegel disk $U$ whose
  rotation number is of bounded type. 
  It is conjectured that the singular value $\kappa$ (and hence
  $\infty$) is accessible both from $U$ and from the escaping set. 
  However, it is currently unknown whether this is true for
  \emph{any} rotation number, even for the golden mean
  (compare Figure \ref{fig:siegel}).

 The analogous problem for quadratic (or, more generally, unicritical)
  polynomials is resolved by showing that the function is quasiconformally
  conjugate to an appropriate model map. In particular, the
  Siegel disk is a quasicircle, and every boundary point is accessible.
  One of the problems
  with producing a similar proof for the exponential family is that 
  it is not clear what the topological model for the boundary $\partial U$
  should be. One might at first expect
  that every component of $\partial U$ is an arc to $\infty$.
  However, the following result shows that the model would need to
  be more complicated.

 \begin{thm}[Indecomposable continua in Siegel disk boundaries] 
   \label{thm:siegeldisks}
  Suppose that $\Ek$ is an exponential map for which the singular
   value $\kappa$ is accessible both from a Siegel disk $U$ and
   from the escaping set $I(\Ek)$. 

  Then $\partial U$ contains uncountably many indecomposable continua.
 \end{thm}

It seems likely that, for every exponential map $\Ek$ with $\kappa\in J(\Ek)$,
 there exists some dynamic ray whose accumulation set contains an entire 
 dynamic ray. Towards this end, we show the following variant
 of Theorem \ref{thm:nonlanding}.

\begin{thm}[Accumulation at Infinity]   \label{thm:accumulation}
 Let $\kappa\in\C$. Then either
  \begin{enumerate}
   \item the accumulation set of every dynamic ray of $\Ek$ is bounded, or
   \item there are uncountably many dynamic rays whose accumulation set
     contains a complete dynamic ray. \label{item:manynonlandingrays}
  \end{enumerate}
\end{thm}
\begin{remark}
 In particular, case (\ref{item:manynonlandingrays}) holds whenever
  the singular value $\kappa$ belongs to the accumulation set of
  some dynamic ray (even if the ray does not land at $\kappa$). 
\end{remark}

\subsection*{Idea of the proof} The proof of the main theorem
  requires a certain amount of combinatorial preparations, so
  let us sketch the underlying idea, which is fairly straightforward.
  Suppose that $g_1$ is some iterated
  preimage of the dynamic ray $g_0$ landing at
  the singular value $\kappa$ of $\Ek$. 
  Since a preimage of a curve landing at $\kappa$ under $\Ek$ will be
  a curve whose real parts tend to $-\infty$, it follows that
  the curve $g_1$ tends to $\infty$ in both directions. 

 Now pick some other dynamic ray $g_2$, very close to $g_1$, which is
  another iterated preimage of $g_0$. Then $g_2$ will stay close to
  $g_1$ for a long time, but eventually tend to $\infty$ in a different 
  direction. If we repeat this process, in the limit we should end up
  with a dynamic ray which does not have a landing point. It is possible
  to ensure that this limit ray will actually accumulate 
  on itself, and to conclude that its closure is an indecomposable
  continuum.

\subsection*{Organization of the article}
 In Section \ref{sec:expcombinatorics}, we review the basic combinatorial
  concepts required for this article, while Section
  \ref{sec:continuity} is devoted to a short discussion of
  the continuity of dynamic rays with respect to their
  external address.
  (A more comprehensive discussion of these topics can
  be found in \cite[Section 2 and 3]{expcombinatorics} and
  \cite[Section 3]{topescapingnew}, respectively.)
  
 In Section \ref{sec:topology}, we give a simple combinatorial
  condition under which the limit set of a 
  ray which accumulates on itself is 
  an indecomposable continuum. 
  Section \ref{sec:nonlanding} contains
  the proof of our main result.

 Finally, Section \ref{sec:accumulation} discusses the accumulation
  of dynamic rays at infinity in a more general setting, in particular
  providing the
  proof of Theorem \ref{thm:accumulation}. Appendix A 
  contains 
  a brief discussion of generalizations of our results
  to larger classes of entire 
  functions.

\subsection*{Acknowledgments}
 The results of this article first appeared as part of
  my thesis \cite{thesis}, and I would like to thank
  my advisor, Walter Bergweiler, for his help and support.
  I would also like to thank
  Adam Epstein, Christian Pommerenke, 
  Johannes R\"uckert, Dierk Schleicher 
  and Sebastian van Strien for
  interesting discussions. 

\subsection*{Notation}
 Throughout this article, $\C$ 
   and $\Ch := \C\cup\{\infty\}$ 
   denote the complex
   plane and the Riemann sphere,
   respectively. The closure of a set
   $A\subset\C$ in $\C$ resp.\
   $\Ch$ will be denoted $\overline{A}$
   and $\widehat{A}$, respectively. 
   The Fatou and Julia sets of an 
    exponential map are denoted
   $F(\Ek)$ and $J(\Ek)$, as usual.    

 Let us also fix the function 
  $F:[0,\infty)\to [0,\infty); t\mapsto \exp(t)-1$
  as a model of exponential growth.  
  We
  conclude any
  proof 
  by the symbol $\blacksquare$, while a result which is cited without
  proof is indicated by $\square$. Separate steps within a proof are 
  concluded by the symbol $\triangle$.

\section{Combinatorics of Exponential Maps} \label{sec:expcombinatorics}
 It has long been customary in exponential dynamics to encode the
  mapping behavior of dynamically defined curves under $\Ek$ using
  symbolic dynamics. We will give a concise summary of these
  concepts here; see \cite[Sections 2 and 3]{expcombinatorics} for 
  a more comprehensive account. 

 A sequence $\s= s_1 s_2 s_3 \dots$ of
  integers is called an \emph{(infinite)
  external address}. 
  If $\Ek$ is an exponential map and $\gamma:(T,\infty)\to\C$ is
  a curve, we say that $\gamma$ \emph{has external address
  $\s$} (as $t\to\infty$) if
  \[ \re\Ek^j(\gamma(t))\underset{t\to\infty}{\to}
     +\infty\quad\text{and}\quad
     \im\Ek^j(\gamma(t))\underset{t\to\infty}{\to} 2\pi s_{j+1} \]
  for all $j\geq 0$. 
  We say that an external
  address $\s$ 
  is \emph{exponentially bounded} if
   \[ \ts := 
        \limsup_{k\to\infty} 
          F^{-(k-1)}(2\pi|s_k|) < \infty. \]
  The space of all exponentially bounded external addresses is denoted
  by $\Sequ_0$.  

 \begin{prop}[Classification of Escaping Points
   {\cite{expescaping}}] \label{prop:rays}
  Let $\kappa\in\C$ and $\s\in\Sequ_0$. Then
   there is a unique maximal curve
   $g_{\s}:(t^{\kappa}_{\s},\infty) \to I(\Ek)$ 
   (where $t^{\kappa}_{\s} \geq t_{\s}$) 
   of escaping points
   which has external address $\s$ as 
   $t\to\infty$ and 
   satisfies 
   \[ |\re \Ek^n(g_{\s}(t)) - F^n(t)| \to 0 \]
   as $n\to \infty$ for any $t>t^{\kappa}_{\s}$.
   This curve
   is unique up to reparametrization 
   and is called the
   \emph{dynamic ray at address $\s$}.
   We say that the ray $\gs$ \emph{lands} at a point
   $z\in \Ch$ if $\lim_{t\to t_{\s}^{\kappa}} \gs(t)=z$. 
 
  If $\kappa\notin I(\Ek)$, 
   then every escaping point is either on a dynamic
   ray or the landing point of a dynamic ray. If $\kappa\in I(\Ek)$, then
   every escaping point eventually maps to a point on a dynamic ray or
   to the landing point of a dynamic ray. \qedd
 \end{prop}

 If $g_{\s}$ is a dynamic ray, we will denote its limit set by
   \[ L_{\s} := 
        \bigcap_{t> t_{\s}^{\kappa}}
           \cl{\gs\bigl([t,\infty)\bigr)}\subset\Ch. \]

 \subsection*{Intermediate external addresses}
  An \emph{intermediate external address}
   is a finite sequence of the form
  \[ \s = s_1 s_2 \dots s_{n-2} s_{n-1} \infty, \]
   where $n\geq 2$, 
   $s_k\in\Z$ for $k<n-1$ and $s_{n-1}\in \Z+ \frac{1}{2}$. 
   The space $\Sequ$ of all infinite and intermediate external addresses,
   equipped with lexicographic order, is order-complete. Its
   one-point compactification is
   $\Sequb := \Sequ\cup \{\infty\}$, which carries
   a complete circular ordering. (We can think of $\infty$ as being
   an intermediate external address of length $1$.) The shift map 
   $\sigma:\Sequ\to\Sequb$ is a locally
   order-preserving map. We will say that $\r_1,\r_2\in\Sequ$
   \emph{surround} an address $\s$ if $\s$
   belongs to the bounded component of
   $\Sequ\setminus\{\r_1,\r_2\}$. 

 \subsection*{Addresses of connected sets}
  Let $\r^-,\r^+\in\Sequ_0$ with $\r^-<\r^+$, and let
   $R>0$. We say that
   $\langle\r^-,\r^+\rangle$ \emph{essentially separates} the half plane
   $\H_R:=\{\re z > R\}$ if the set
   \[ \H_R\setminus \bigl(
           g_{\r^-}\bigl([t_{\r^-}^{\kappa}+1,\infty)\bigr) \cup
           g_{\r^+}\bigl([t_{\r^+}^{\kappa}+1,\infty)\bigr) \bigr) \]
  has a component $U$ with unbounded real parts but bounded imaginary
   parts. (In other words, both ray pieces intersect the
   line $\{\re z = R\}$.)
   The component $U$ (if it exists) is necessarily unique, and 
   will be denoted by $U_R(\langle\r^-,\r^+\rangle)$. 

  Let $A\subset\C$, and let
   $\s\in\Sequ_0$. We say that 
   \emph{$\s$ is separated from
   $A$} if there exist
   $\r^-,\r^+\in \Sequ_0$ and some 
   $R>0$ such that $\r^-$ and $\r^+$ surround $\s$,
   the pair $\langle\r^-,\r^+\rangle$ essentially
   separates $\H_R$ and
   $A\cap U_R(\langle\r^-,\r^+\rangle) = \emptyset$. Similarly, we say that
   the address
   $\infty\in\Sequb$ is separated from $A$ if $\im A$ is bounded and
   $\re A$ is bounded from below. 
   The set
   \[ \ADDR(A) := 
        \{\s\in\Sequb: \text{$A$ is not separated from $\s$}\} \]
   is clearly compact, and is empty if and only if
   $A$ is bounded. 

   \begin{remark}
    Another way to phrase this definition is as follows: we form
     a natural compactification of $\C$ by adjoining the space
     $\Sequb$ as a circle at $\infty$. The set
     $\ADDR(A)$ is then exactly the set of accumulation points of $A$ in
     $\Sequb$ with respect to this topology.
   \end{remark}

 \subsection*{Two directions of dynamic rays}  
  We will usually apply the above concepts to
  dynamic rays $\gs$ which accumulate
  at
  $\infty$ as $t\searrow t_{\s}^{\kappa}$. In order to 
  facilitate these discussions, let us abbreviate
   \[ \ADDR^-(\s) := 
      \ADDR^-(\gs) := 
      \ADDR\bigl(\gs\bigl((t_{\s}^{\kappa},t_{\s}^{\kappa}+1]\bigr) \bigr). \]

 \subsection*{Itineraries}

 \begin{defn}[Accessible Singular Values]
  Let $\Ek$ be an exponential map with $\kappa\in J(\Ek)$. We say that
  the singular value is \emph{accessible} if $\kappa$ is either on a 
  dynamic
  ray or the landing point of a dynamic ray. If $\r$ is the external address of
  such a ray, we write $\r = \extaddr(\kappa)$.
 \end{defn}
 \begin{remark}
 \begin{enumerate}
  \item[1.]
   The address $\r$ need not be unique. 
  \item[2.]
   Misiurewicz and escaping parameters have accessible singular values, as
   mentioned in the introduction. Furthermore, the landing points of
   parameter rays at ``regular''
   addresses in the sense of Devaney, Goldberg and Hubbard \cite{dghnew1}
   have this property. We expect that this condition holds
   in most situations in which
   local connectivity is known for unicritical polynomials, for example 
  the Siegel parameters mentioned in the introduction.
 In particular, we believe that accessibility
   should be generic in the bifurcation locus and includes many parameters
   for which the singular orbit is dense in the plane.
  \end{enumerate}
 \end{remark}

 Suppose that $\Ek$ is an exponential map with accessible singular value
  in $J(\Ek)$,  
  and let $\r = \extaddr(\kappa)$. Then the preimages
  $\Ek^{-1}(\gr)$ (that is, the dynamic rays at addresses of the form
  $k\r$ for $k\in\Z$) cut the plane into countably many strips $S_k$. 
  Let us label these such that $S_k$ is the strip bounded by
  $g_{k\r}$ and $g_{(k+1)\r}$. (We refer to these strips as the
  \emph{dynamic partition}.)

 If $z\in\C$, we can assign to $z$ an \emph{itinerary}
   \[ \itin(z) := \itin_{\r}(z) := \addu := \u_1 \u_2 \u_3 \dots \]
 where $\u_j = k$ if $\Ek^{j-1}(z)\in S_k$, and 
   $\u_j = \left(\itj\right)$
 if $\Ek^{j-1}(z)\in g_{j\r}$. 

 If $z\in g_{\s}$ for some $\s\in\Sequ_0$, then the itinerary entries
  of $z$ clearly satisfy
   \begin{equation}
    \u_k = \begin{cases}
              \j & \text{if } \j \r < \sigma^{k-1}(\s) < (\j+1) \r \\
              \left(\itj\right) & \text{if }
                                            \sigma^{k-1}(\s) = \j \r.
           \end{cases} \label{eqn:itin}
   \end{equation}
 
 For any $\r\in\Sequ$ and every infinite external address $\s$, we can
  define an address
  $\itin_{\r}(\s)$ by the formula (\ref{eqn:itin}). If
  $\s$ is an intermediate external address of length $n$, then we similarly
  define
  \[ \itin_{\r}(\s) = \u_1 \dots \u_{n-1} *, \]
  where 
  $\u_1,\dots,\u_{n-1}\in\Z$ satisfy (\ref{eqn:itin}). 
  The \emph{kneading sequence} of $\s\in\Sequ$ is
  the itinerary $\K(\s) := \itin_{\s}(\s)$

  Finally, let us say that an itinerary entry 
  $\m\in\Z$ is \emph{adjacent} to the entry $\u$ if
  $\m=\u$ or $\u$ is a boundary sumbol $\left(\itj\right)$ with
  $\m\in\{\j-1,\j\}$.

 We will frequently use the following simple observation. 

\begin{lem}[Addresses sharing an itinerary] \label{lem:itineraries}
 Let $\s\in\Sequ$ and
  $\addu := \K(\s)$.  Suppose that
  $\r\neq \wt{\r}$ are two addresses
  sharing the same itinerary $\addut := \itin_{\s}(\r)=\itin_{\s}(\wt{\r})$,
  and let $m\geq 1$ with $r_m\neq \wt{r}_m$. 
  
 Then for  
  every $k\geq 0$, there exists $0\leq j\leq k$ such that 
  $\sigma^{m+k}(\r)$ and $\sigma^{m+k}(\wt{\r})$ surround
  $\sigma^j(\s)$. 
  In particular, $\ut_{m+k}\in \{\u_1,\dots,\u_k\}$ for all
  $k\geq 1$. 
\end{lem}
\begin{proof} Let $k\geq 0$. 
  We may 
  suppose without loss of generality 
  that
  $r_{m+i}=\wt{r}_{m+i}$ for all 
  $i\in\{1,\dots,k\}$.
  (Otherwise, we can replace
   $m$ by $m+i$ and $k$ by $k-i$.)

  By the definition of itineraries, $\s$ belongs to 
  the interval $I$ of $\Sequ$ between
  $\sigma^m(\r)$ and $\sigma^m(\wt{\r})$. Since we assumed that 
  the latter addresses
  agree in their first $k$ entries, it follows that
  $\sigma^k(I)\ni \sigma^k(\s)$ is the interval bounded by
  $\sigma^{m+k}(\r)$ and $\sigma^{m+k}(\wt{\r})$, as required. \end{proof}

\section{Continuity among Dynamic Rays} \label{sec:continuity}
 The dependence of $\gs(t)$ on the pair $(\s,t)$ is quite complicated
  in general, and will depend to some degree on the parameter
  $\kappa$. On the other hand, it is well-known that there is 
  a certain continuity in the $\s$-direction, e.g.\
  in the sense that every
  compact subpiece of the dynamic ray $\gs$ can be approximated
  from above and below by suitable sequences of dynamic rays. In fact,
  in \cite{topescapingnew}, a \emph{simultaneous} parametrization of
  all dynamic rays was given, which greatly simplifies both questions
  of this kind and those regarding  escaping 
  endpoints of rays. However, to reduce the prerequisites for this
  article, and to improve consistency among
  articles on exponential dynamics, we shall state the required
  results here using the
  original parametrization from \cite{expescaping}, and sketch a
  proof using only results from \cite{expescaping}.

  \begin{prop}[Asymptotics of Dynamic Rays 
                  {\cite[Proposition 3.4]{expescaping}}]
    \label{prop:asymptotics}
   Let $\kappa\in\C$ and $\s\in\Sequ_0$, and set
  \[  \ts^* := \sup_{k\in\N} F^{-(k-1)}(2\pi|s_k|). \]
   Then $\ts^{\kappa} < T_{\s} := 2\ts^* + \log^+|\kappa|+4$, and 
    for all $t\geq T_{\s}$,    
    \[ \bigl|\gs(t) - (t + 2\pi i s_1)\bigr| < 
           e^{-t/2}. \qedd \]
  \end{prop}
  \begin{remark}
   This is not quite the formulation given in
    \cite{expescaping}, but is easily deduced from it.
  \end{remark}

 \begin{lem}[Continuity Between Rays] \label{lem:raycontinuity}
  Let $\kappa\in\C$,
   $\s\in\Sequ_0$ and $K>0$. For $n_0\in\N$,
   denote by $S(\s,K,n_0)$ the
   set of all addresses $\wt{\s}$ which 
   agree with $\s$ in the first $n_0$ entries and satisfy
    $|\wt{s}_n| \leq K + \max_{k\leq n}|s_k|$
    for all $n\geq 1$. 

   Then for every $t_0>\ts^{\kappa}$ and every $\eps>0$, there exists
    $n_0$ such that $t_0 > t_{\tilde{\s}}^{\kappa}$ and 
    \[ |\gs(t) - g_{\tilde{\s}}(t)| < \eps \]
    for all $t\geq t_0$ and all $\wt{\s}\in S(\s,K,n_0)$. 
 \end{lem}
 \begin{fullsketch}
  It is easy to see (from the definitions of
   $\ts$ and $\ts^*$, together with
   the fact that $t_0 > \ts^{\kappa}\geq \ts$),
   that there is some $n_1$ such that
   \[ F^n(t_0) > 2t^*_{\sigma^n(\s)}+2K + |\kappa|+4 \]
   for all $n\geq n_1$. By definition, we have 
   $t^*_{\sigma^n(\tilde{\s})} \leq t^*_{\sigma^n(\s)}+K$
   for all $\s\in S(\s,K,n)$. 
   Thus we have
   $F^n(t) > T_{\sigma^n(\tilde{\s})}$ (where
   $T_{\sigma^n(\tilde{\s})}$ is
    the number from Proposition \ref{prop:asymptotics})
   whenever $n\geq n_1$,
   $\wt{\s}\in S(\s,K,n+1)$ and $t\geq t_0$.
   In particular, since $F^n(t)>|\kappa|+4$, we have
   \begin{align*}
       &\re g_{\sigma^n(\tilde{\s})}\bigl(F^n(t)\bigr) 
       \geq F^n(t)-e^{-F^n(t)/2} > \re\kappa + 2 \quad\text{and} \\
       |&g_{\sigma^n(\tilde{\s})}(F^n(t)) - g_{\sigma^n(\s)}(F^n(t))| < 
        2e^{-F^n(t)/2} < 1.
   \end{align*}

  In particular, for $n\geq n_1$, the pieces 
   $g^n := g_{\sigma^n(\s)}\bigl([F^n(t_0),\infty)\bigr)$ are contained
   in the half-plane $H := \{z: \re(z-\kappa)> 2\}$. Note that
   $\Ek$ is expanding on 
   $H$. It follows easily that, 
   for sufficiently large $n$, there exists a branch
    \[ \phi: \bigl\{z\in\C: \dist(z,g^n)<1 \bigr\} \to \C \]
   of $\Ek^{-n}$ with $\phi(g_{\sigma^n(\s)}(F^n(t)))=g_{\s}(t)$ for
   $t\geq t_0$ such that
   $|\phi'| < \eps$. By the definition of 
   dynamic rays, we have 
   \[ \phi(g_{\sigma^n(\tilde{\s})}(F^n(t))) =
        g_{\tilde{\s}}(t) \]
   for all $\wt{\s}\in S(\s,K,n+1)$,
   and the claims follow. 
 \end{fullsketch}

   \begin{lem}[Accumulation on 
                      Dynamic Rays] 
              \label{lem:accumulation} 
    Let $\Ek$ be
    an exponential map, and
    suppose that $A\subset\C$ is 
    connected and intersects at most 
     countably many dynamic rays
     of $\Ek$. If $\s\in \Sequ_0$ with 
     $\s\in \ADDR(A)$, then $\gs\subset \cl{A}$ or 
     $\cl{A}\subset \gs$. 
   \end{lem}
   \begin{proof}
   Since $A$ only intersects countably 
   many dynamic rays, we can find
   a sequence $(k_n)$ with $k_n\in\{1,2\}$ such that, for all $n$,
   $A$ does not intersect the dynamic rays
    $g_{\r_n^+}$ and $g_{\r_n^-}$ defined by
   \[ \r_n^{\pm} := s_1 \dots s_n (s_{n+1} \pm k_{n+1})
                 (s_{n+2} \pm k_{n+1}) \dots. \]
   By Lemma \ref{lem:raycontinuity},
   $g_{\r_n^+} \to g_{\s}$
   uniformly on every interval $[t,\infty)$ with $t>
   \ts^{\kappa}$, and the same is true for $\r_n^-$. 

   So the rays $g_{\r_n^{\pm}}$ approximate $\gs$ from above and below.
    Suppose that there is
    $t_0>\ts^{\kappa}$ such
    that $\eps := 
              \dist(\gs(t_0),A)>0$,
    and let $U_n$ denote the
    component of
    \[ \C\setminus \bigl( 
       \D_{\eps}(\gs(t_0))
        \cup g_{\r_n^+}\cup
        g_{\r_n^-} \bigr) \]
     which contains $\gs(t)$ for
     large $t$. Then $A\subset U_n$.
     Clearly 
     $\bigcap U_n \subset 
      \gs\bigl([t_0,\infty)\bigr)$,
     and the claim follows.
   \end{proof}

  \begin{cor}   \label{cor:addresses}
   Let $\Ek$ have an accessible singular value, and
    $\s=\extaddr(\kappa)$, and let $\r\in\Sequ_0$. 
   \begin{enumerate}
    \item Suppose that
       $\Addr^-(\r)$ contains an intermediate external address of
       length $m$ 
       or a preimage of $\s$ under $\sigma^m$. 
       Then the ray $g_{\sigma^m(\r)}$ accumulates on $\kappa$. 
       In particular, all itinerary entries of $\sigma^m(\r)$ are
       adjacent to the corresponding entries of $\K(\s)$. 
    \item If the ray $\gr$ lands at $\infty$, then 
      $g_{\sigma^m(\r)}$ lands at $\kappa$ for some $m > 0$. 
   \end{enumerate}
  \end{cor}
 \begin{proof} If $\Addr^-(\r)$ contains
      a preimage of $\s$ under $\sigma^m$, then
      $\s\in\Addr^-(\sigma^m(\r))$. By the previous lemma,
      $L_{\sigma^m(\r)}$ contains the entire ray
      $g_{\s}$, and therefore also $\kappa$, its landing point. 
     
    Similarly, if $\Addr^-(\r)$ contains an intermediate external
     address of length $m$, then $\infty\in \Addr^-(\sigma^{m-1}(\r))$. 
     Since $g_{\sigma^{m-1}(\r)}$ belongs to some strip of the dynamic
     partition, this is only possible if $\im(g_{\sigma^{m-1}(\r)})$ 
     is unbounded from below, which means that
     $g_{\sigma^m(\r)}$ accumulates at $\kappa$. 

     Finally, suppose that $g_{\r}$ lands at $\infty$. Then 
      $\Addr^-(\r)$ consists of a single address $\wt{\r}$. By the
      previous lemma, $\wt{\r}\notin\Sequ_0$, as otherwise the accumulation
      set $L_{\r}$ would contain an entire dynamic ray. Also,
      $\wt{\r}$ cannot be an exponentially unbounded infinite address
      since every entry of its
      itinerary must be adjacent to the corresponding
      entry of $\itin_{\s}(\r)$. So
      $\wt{\r}$ is an intermediate external address of length $m$, 
      which means that 
      $g_{\sigma^m(\r)}$ lands at $\kappa$. \end{proof}

\section{Topological Considerations} \label{sec:topology}

A useful tool for showing that the accumulation set of a dynamic ray
is an indecomposable continuum is given by the following theorem of
Curry \cite{curry}.

\begin{thm}[Curry] \label{thm:curry}
 Suppose that $g$ is a curve in $\Ch$, and let $G$ denote its
 accumulation set. If $G$ has
 topological dimension one, does not separate the Riemann
 sphere into infinitely many components
 and contains $g$, then $G$ is an indecomposable continuum. \qedd
\end{thm}

\begin{lem}[Accumulation sets of dynamic rays] \label{lem:indecomposable}
 Let $\Ek$ be an exponential map with accessible singular value. 
  Then for every $\r\in\Sequ_0$, the
  accumulation set $L_{\r}$ has empty interior.

 Furthermore, suppose that $\Addr^-(\r)$ is finite. 
  Then $\C\setminus L_{\r}$ has only finitely many
  components. 
\end{lem}
\begin{proof} 
 Let  $\s=\extaddr(\kappa)$ and
  $\ut:= \itin_{\s}(\r)$. Then
  for every $n\geq 0$, the set 
  $\Ek^n(L_{\r})$ is contained in the closure of some 
  strip
  $S_{\ut_{n+1}}$ of the dynamic partition. Thus we can find some
  $\m\in\Z$ and a subsequence $n_j$ such that 
  $\Ek^{n_j}(\gs)\cap S_{\m}=\emptyset$ for all $j$. If $U$ was a
  component of $\operatorname{int}(L_{\r})$, 
  then $\Ek^{n_j}|_U$ would omit all points of
  $S_{\m}$, and thus $U\subset F(\Ek)$ by Montel's theorem. 
  This contradicts the fact that $L_{\r}\subset J(\Ek)$. 

 Now suppose that $\Addr^-(\r)$ is finite. Then for each component
  $I$ of $\Sequ\setminus \Addr^-(\r)$, there is a unique component $V$
  of $\C\setminus L_{\r}$ such that, for every
  $\wt{\r}\in I$, the ray $g_{\wt{\r}}$ eventually tends to $\infty$ in $V$.
  Since $\Sequ\setminus \Addr^-(\r)$ is finite, 
  there are only finitely many such
  components $V$. Let $U$ denote the union of all other components of
  $\C\setminus L_{\r}$; we claim that $U$ has at most one component. 

 Similarly as above, it follows that $U\subset F(f)$. 
  Since $\kappa\in J(f)$, it follows that every component of $U$ is 
  a preimage of a Siegel disk $V$ of $\Ek$. By passing to a forward
  iterate, we may suppose without loss of generality that $V\subset U$. 
  Since $\Ek$ is injective on $U$, it is impossible for $U$ to
  contain any other component which eventually maps to $V$. \end{proof}

\section{Existence of Nonlanding Rays} \label{sec:nonlanding}
 \begin{thm} \label{thm:nonlandingprecise}
  Let $\Ek$ be an exponential map 
  with accessible singular value
  $\kappa\in J(\Ek)$. Then there exists an uncountable
  set $R\subset \Sequ_0$ such that
  \begin{enumerate}
   \item $\Addr^-(\r)=\{\r\}$ for all $\r\in R$, and
     \label{item:accumulation}   
   \item no two addresses in $R$ share the same itinerary.
     \label{item:differentitineraries}
  \end{enumerate}
  If $\extaddr(\kappa)$ is bounded, then $R$ can be chosen to consist only
   of bounded addresses.
 \end{thm}
 \begin{proof}
  Let $\s := \extaddr(\kappa)$ and set
   \[ T_n := 2 + \max_{k\leq n} s_k \]
  for every $n\in\N$. 
  We define $R_1$ to be the set of all addresses of the form
  \begin{align} \label{eqn:nonlandingaddress}
    \r &=  \r(n_1,n_2,n_3,\dots) \\
       & := T_1 
         s_1 s_2 \dots s_{n_1-1} T_{n_1} 
         s_1 s_2 \dots s_{n_2-1} T_{n_2} 
         s_1 s_2 \dots
         s_{n_3-1}T_{n_3} \dots \ , \notag
  \end{align}
  where $(n_k)$ is some sequence of natural numbers. 
  (The set $R$ will be a suitable
   subset of
   $R_1$.)
  Given 
  $N_1,\dots,N_k\in\N$, let us also denote by 
  $R_1(N_1,\dots,N_k)$ the subset of $R_1$ consisting of all sequences
  $\r(n_1,n_2,\dots)$ with $n_j=N_j$ for $j\leq k$.

  \begin{claim}[Claim 1] For every $\r\in R_1$, 
     $t_{\r}^*\leq t_{\s}^*+2$. In particular,
     $t_{\r}^{\kappa}\leq T_0 := 2t_{\s}^*+\log^+|\kappa|+8$.
  \end{claim}
 \begin{subproof}
  This follows from the definitions and 
   Proposition \ref{prop:asymptotics}. 
 \end{subproof}

 \begin{claim}[Claim 2]
   For every $\r\in R_1$, 
    there are no other addresses whose itinerary
    coincides with that of $\r$. In particular,
    no two addresses in $R_1$ share the same itinerary and
    $\Addr^-(\r)\subset \{\r\}$ for all $\r\in R_1$. 
  \end{claim}
  \begin{subproof}
   Set $\addu := \K(\s)$ and
    $\addut := \itin_{\s}(\r)$. By definition of
    $R_1$, for every $m\in\N$ there is some $k\geq 1$ such that
    $r_{m+k} = T_{k'}$ for some $k'\geq k$. In particular, 
    $r_{m+k} > 2+s_{k}$, and therefore 
    $\u_k \neq \ut_{m+k}$. By Lemma \ref{lem:itineraries}, there can
    be no other address with the same itinerary.
    
   In particular, $\Addr^-(\r)$ cannot contain any infinite
    external addresses other
    than $\r$ which are not preimages of $\s$. On the other hand,
    $\Addr^-(\r)$ also cannot contain any preimages of $\s$ or
    any intermediate external addresses by Corollary \ref{cor:addresses}.
  \end{subproof}

  We define $R$ to be the subset of all $\r\in R_1$ with
   $\Addr^-(\r)\neq\emptyset$. We will show that $R$ is nonempty
   by inductively constructing a sequence $(n_k)$ for which
   the dynamic ray at address $\r(n_1,n_2,\dots)$ accumulates at
   infinity. In this construction, there will 
   be countably many different choices for
   each $n_j$, so that $R$ is in fact uncountable, as claimed.

  Suppose that $j\geq 1$ such that $n_k$ has been chosen for all
   $k<j$. Let $\r^k$ denote the address
   \[ \r^j :=
       \r(n_1,\dots,n_{j-1}) := T_1 s_1 s_2 \dots
                  s_{n_1-1}T_{n_1} s_1\dots
                  s_{n_{j-1}-1}  T_{n_{j-1}} \s. \]
   Then $\r^j$ is a preimage of $\s$, and thus lands at
    $\infty$. Let
    $t_j\in (t_{\r^j}, T_0]$ be the largest value satisfying 
    $|g_{\r^j}(t_j)| > j$. By Lemma
    \ref{lem:raycontinuity}, there is some $N_j\in\N$ such that
    \[ |g_{\r}(t_j)|\geq j \]
    for all $\r\in R_1(n_1,\dots,n_{j-1},n)$ with
    $n \geq N_j$. We choose $n_j$ to be any such $n$. 

   Let $\r=\r(n_1,n_2,\dots)$ be an address constructed in this way. 
    Then
    $g_{\r}(t_j)\to \infty$.
    Since $t_j \leq T_0$ for all $j$, we must have $t_j\to t_{\r}^{\kappa}$.
    I.e., $g_{\r}$ accumulates at $\infty$ and thus $\r\in R$,
    as required. \end{proof}

\begin{proof}[Proof of Theorem \ref{thm:nonlanding}]
  This theorem is a direct corollary of the previous theorem
   together with Lemma \ref{lem:indecomposable} and
   Theorem \ref{thm:curry}.
\end{proof}

\begin{remark}
 Our construction yields nonseparating
  indecomposable continua, each containing
  exactly one dynamic ray. For many addresses $\addr(\kappa)$, it is 
  possible to modify the above construction to obtain two addresses
  $\r^1,\r^2$ with $\ADDR^-(\r^1)=\ADDR^-(\r^2)=\{\r^1,\r^2\}$. This
  leads to an indecomposable limit set which separates the plane into
  two components and contains the two dynamic rays $g_{\r^1}$ and
  $g_{\r^2}$; see also 
  \cite{misindecomposable}, where this is carried out for
  the case of the Misiurewicz parameter 
  $\kappa=\log(2\pi)+\pi i/2$ with $\addr(\kappa)=0111\dots$.
 (For this particular parameter,
  \cite{misindecomposable} also constructs addresses
  $\r^1$ and $\r^2$ with $\ADDR^-(\r^1)=\ADDR^-(\r^2)=\{\r^2\}$, so that
  $g_{\r^1}$ accumulates on an indecomposable continuum but not on
  itself.)
\end{remark}

\begin{proof}[Proof of Theorem \ref{thm:nopincheddisk}]
 For attracting or parabolic
 parameters, every dynamic ray lands and
 the Julia set is the union of these rays and their landing points
 \cite[Theorem 5.7]{accessible} (compare also
 \cite[Corollary 9.3]{topescapingnew} for the stronger statement that
 the dynamics on $J(\Ek)$ is semiconjugate to that of an
 exponential map with an attracting fixed point). 
 For the converse direction, let $\Ek$ be an 
 exponential map with $\kappa\in J(\Ek)$. Then one of the following
 holds:
  \begin{enumerate}
   \item The singular value $\kappa$ is not accessible from
     $I(\Ek)$; i.e., $\kappa$ is neither on a dynamic
     ray nor the endpoint of
     a ray, or
   \item The singular value $\kappa$ is accessible from
     $I(\Ek)$, in which case there exists a nonlanding 
     dynamic ray of $\Ek$ by Theorem \ref{thm:nonlanding}.\qedhere
  \end{enumerate}
\end{proof} 

\begin{proof}[Proof of Corollary \ref{cor:uncountablymanydisjointcontinua}]
 Since the indecomposable continua constructed in Theorem
  \ref{thm:nonlandingprecise} all have different itineraries, they only
  have the point at $\infty$ in common. As indicated in Section
  \ref{sec:expcombinatorics}, 
  it is not difficult to compactify $\C$ to a space 
  $\wt{\C}$ in such a way that each external address $\s\in\Sequb$ 
  corresponds to (exactly) one point at $\infty$. So if, for $\r\in R$,
  we take the closures
  $\wt{g_{\r}}$ in $\wt{\C}$, the resulting indecomposable continua 
  are pairwise disjoint. The space $\wt{\C}$ is clearly homeomorphic to
  the closed unit disk $\cl{\D}\subset\C$, concluding the proof.
\end{proof}

 Let us conclude the section by indicating how the above construction can be
  modified to yield a proof of Theorem \ref{thm:siegeldisks}
  (concerning Siegel disks). Suppose that $\Ek$ has a Siegel disk $U$ whose
  boundary contains the singular value, and suppose furthermore that
  $\kappa$ is accessible from the escaping set; say $\extaddr(\kappa)=\s$.
  Then clearly $\u := \K(\s)$ is periodic, but 
  $\s$ is not. The period of $\u$ is at most the period of the Siegel
  disk $U$. It is less obvious, but also true, that these periods must be
  equal. (This follows, for example, from the combinatorial results of
  \cite{expcombinatorics}; we omit the details here.) 
 \begin{thm}
   Let $\s\in\Sequ_0$ be a non-periodic external address for which
     $\addu := \K(\s)$ is periodic. Suppose that $\kappa$ is a parameter with 
     $\extaddr(\kappa)=\s$. Then there exist uncountably many addresses
     $\r\in\Sequ_0$ with $\itin_{\s}(\r)=\addu$ and 
     $\r\in\ADDR^-(\r)$.

   Suppose furthermore that $\Ek$ has an unbounded Siegel disk $U$. Then,
    for each of these addresses, 
    $g_{\r}\subset \partial U$. If $\kappa$ is accessible from $U$, then
    $\ADDR^-(\r)=\{\r\}$; in particular, $L_{\r}$ is an indecomposable
    continuum. 
 \end{thm}
 \begin{fullsketch}
 For simplicity, we will restrict to the case where
  $U$ is a fixed Siegel disk (the general case is analogous, but requires 
  slightly more bookkeeping). Shifting $\kappa$
   by an integer multiple of $2\pi i$, we may then furthermore assume that
   $\addu=\texttt{000}\dots$. 

  Let $X\subset\Sequ_0$ be the compact set
   consisting of all infinite external addresses $\r$
   whose itineraries are adjacent to $\addu$; i.e.\
   $\addu\in\{\itin^+(\r),\itin^-(\r)\}$.
    We claim that every interval of $\Sequb\setminus X$ is
    bounded by two iterated preimages of $\s$ under the shift, and
    contains an intermediate external address $\t$ with
    $\itin(\t)=\texttt{00}\dots\texttt{00}*$. Indeed, this follows easily
    from Lemma \ref{lem:itineraries} and the fact that $\s\in X$. 
    In particular,
    iterated preimages of $\s$ are dense in $X$. 
    (In fact, the map which collapses the two
     addresses $0\s$ and $1\s$, and every pair of their preimages,
     semi-conjugates $\sigma|_X$  to
     an irrational rotation of the circle. In particular, \emph{every}
     backward (and forward)
     orbit is dense in $X$.)

   Now we can inductively construct uncountably many
    addresses $\r\in X$ with 
    $\Addr^-(\r)\neq\emptyset$: 
    as in the proof of Theorem \ref{thm:nonlandingprecise},  
    $\r$ is the limit of
    iterated preimages $\r^j\in X$ 
    of $\s$, where $\r^{j+1}$ is chosen sufficiently
    close to $\r^j$.

   To obtain the stronger claim, stated in the theorem, that
    actually $\r\in\Addr^-(\r)$, let $\r^{j-}$ be the intermediate
    address which is the unique
    element of $\Addr^-(\r^j)$. 
    By choosing $\r^{j+1}$ sufficiently close to
    $\r^j$ in each step, we can easily ensure that $\Addr^-(\r)$
    contains a limit address of the sequence $\r^{j-}$.
    However, 
    it is not difficult to show
    (using Lemma \ref{lem:itineraries}) that $\r^{j-}\to \r$.
     Hence $\r\in\Addr^-(\r)$, which completes the
    proof of the first claim.

   Now suppose that $\Ek$ has an unbounded, fixed Siegel disk $U$,
    and consider the set $\Addr(U)$. By \cite{siegel}, we have
    $\kappa\in\partial U$. In particular, $U$ is unbounded to the left, so
    $\infty\in\Addr(U)$. It follows that any intermediate external address
     $\t$ with
     $\itin(\t)=\texttt{00}\dots\texttt{00}*$ also belongs to 
     $\Addr(U)$. By the above observation (and since $\Addr(U)$ is compact),
     this implies that $X\subset\Addr(U)$, and hence
     $\gr\subset\cl{U}$ for every $\r\in X$. 

    Finally, suppose that
     $\kappa$ is accessible from $U$, say by a curve $\gamma$
     connecting it to the indifferent fixed point. Then, for any
     $\r\in X$ which is not an iterated preimage of $\s$, we can use
     preimages of $\gamma$ to separate $\gr$ from all other dynamic rays.
     This shows that $\Addr^-(\r)\subset\{\r\}$, as desired.
    If $\infty$ is accessible from $U$ by a curve $\gamma$, but $\kappa$
     is not, then $\t := \extaddr(\gamma)\in X$. The iterated preimages of
     $\t$ are dense in $X$, and the claim follows analogously.
 \end{fullsketch}

\section{Exponential maps with rays accumulating at infinity}
  \label{sec:accumulation}

 In the previous sections, we dealt with the case where
  $\kappa$ is accessible. In this section, we will briefly
  explore the situation where $\kappa$ is merely an accumulation
  point of some dynamic ray, or even more generally
  when we merely know that some
  dynamic ray accumulates at $\infty$. The latter case is equivalent
  to saying that there is some $\s\in\Sequ_0$ with $\Addr^-(\s)\neq\emptyset$;
  let us begin by showing that there is a similar criterion for the
  former case.

 \begin{lem}[Rays accumulating at the singular value]
  Let $\kappa\in\C$ and $\s\in\Sequ_0$.
   Then $g_{\sigma(\s)}$ accumulates on the singular value
   $\kappa$ if and only if 
   $\infty\in\Addr^-(\s)$. 
 \end{lem}
 \begin{proof} (Compare also
  \cite{siegel}.)
  We may assume that $\kappa\in J(\Ek)$, since otherwise
  all dynamic rays land in $\C$. The ``only if'' part of the
  statement is trivial. So
  let $\s\in\Sequ_0$ such that $g_{\sigma(\s)}$ does not
  accumulate on $\kappa$, and let $U$ be the component of
  $\C\setminus L_{\sigma(\s)}$ containing $\kappa$. 
  Since $\kappa\in J(f)$, the set $U$ contains some 
  escaping point $z_0$, say $z_0=g_{\r}(t_0)$. Consequently
  $\gamma_0 := g_{\r}\bigl([t_0,\infty)\bigr)\subset U$. Extend the
  curve $\gamma_0$ to a curve $\gamma\subset U$ by connecting 
  $\kappa$ and $z_0$. Then there is a branch of $\Ek^{-1}$ defined on
  $\C\setminus \gamma$ taking $g_{\sigma(\s)}$ to $g_{\s}$. It follows
  easily that $\Addr^-(\s)$ is contained in the
  interval $(m\r,(m+1)\r)$ of $\Sequ$ for some $m\in\Z$;
  in particular, $\infty\notin\Addr^-(\s)$. \end{proof}
 
 Let us now investigate what happens when we attempt to adapt the
  proof of Theorem \ref{thm:nonlandingprecise} to our more
  general situation. Again, there are two cases to consider.
  If there is some address $\s\in\Sequ_0$ for which
  $\kappa\in L_{\s}$ and $L_{\s}$ is bounded, then we still have
  the tool of itineraries at our disposal, and the combinatorial
  part of the proof will go through just like before. Otherwise,
  we will still be able to construct many rays $g_{\r}$ 
  which accumulate
  at infinity, but without much control over $\Addr^-(\r)$. 

 \begin{thm}[Rays accumulating at $\infty$] \label{thm:raysatinfty}
  Let $\kappa\in\C$. 
  \begin{enumerate}
   \item The set of addresses $\r\in\Sequ_0$ with
     $\Addr^-(\r)\neq\emptyset$ is either empty or uncountable.
       \label{item:uncountablymanyaddresses}
   \item 
    If there is an address $\s\in\Sequ_0$ with
     $\Addr^-(\s)=\emptyset$ such that $\kappa\in L_{\s}$,
     then there are
     uncountably many addresses $\r\in\Sequ_0$ with
     $\Addr^-(\r)=\{\r\}$. \label{item:addressinfty}
  \end{enumerate}
 \end{thm}
 \begin{proof} 
   If there is an address $\s\in\Sequ_0$ with
    $\Addr^-(\s)\neq\emptyset$, then we can construct
    uncountably many addresses with $\Addr^-(\r)\neq\emptyset$
    as in (the second part of) the proof of
    Theorem \ref{thm:nonlandingprecise}, as a limit of preimages of
    $\s$. If furthermore (\ref{item:addressinfty}) holds, then 
    the first part of the proof of Theorem \ref{thm:nonlandingprecise}
    will also go through, allowing us to construct addresses with
    $\Addr^-(\r)=\{\r\}$.
  \end{proof}

  It is not difficult
   to modify the proof of (\ref{item:uncountablymanyaddresses}) 
   to directly yield the existence of
   uncountably many addresses $\r$ for which
   $\Addr^-(\r)\cap \Sequ_0\neq \emptyset$, proving
   Theorem \ref{thm:accumulation}. Instead, we will 
   use a little more combinatorics to obtain a slightly stronger
   statement: if there is no address $\s$ as in
   Theorem \ref{thm:raysatinfty} (\ref{item:addressinfty}),
   then every $\r\in\Sequ_0$ with $\Addr^-(\r)\neq\emptyset$
   satisfies $\Addr^-(\r)\cap \Sequ_0\neq\emptyset$.

 \begin{lem}[Close accumulation addresses]
  Let $\kappa\in\C$, and let $\s\in\Sequ_0$. Then
   for every $\r\in\Addr^-(\sigma(\s))$, there is 
    $\wt{\r}\in\Addr^-(\s)$ with $\sigma(\wt{\r})=\r$ and
    $|\wt{r}_1 - s_1|\leq 1$. 
  
  In particular,
   $\Addr^-(\sigma(\s))=\sigma(\Addr^-(\s)\setminus\{\infty\})$.
 \end{lem}
 \begin{proof} The first statement is a consequence of the fact that
  $g_{\s}$ cannot intersect its $2\pi i\Z$-translates. 
  To prove the second statement, note that 
  $\Addr^-(\sigma(\s))\subset \sigma(\Addr^-(\s)\setminus\{\infty\})$
  follows from the first claim, and
  the converse inclusion follows from the definitions. \end{proof}

 \begin{cor} \label{cor:boundeddistance}
  Let $\Ek$ be an exponential map and let $\s\in \Sequ_0$ and
   $m\geq 1$ such that
   $\Addr^-(\sigma^m(\s))\neq \emptyset$. Then 
   there is $\wt{\s}\in\Addr^-(\s)$ satisfying
    $|s_j - \wt{s}_j|\leq 1$ for $j\in\{1,\dots,m\}$. 
 \end{cor}
 \begin{proof} This follows from the previous lemma by a simple
  induction. \end{proof}

 \begin{lem}[Exponentially bounded accumulation addresses]
  Let $\Ek$ be an exponential map and let $\s\in\Sequ_0$. Then one
   of the following holds:
  \begin{enumerate}
   \item $\Addr^-(\s)\cap \Sequ_0\neq \emptyset$, or
   \item $\Addr^-(\sigma^m(\s))=\emptyset$ for some $m\geq 0$.
  \end{enumerate}
 \end{lem}
 \begin{proof} Suppose that $\Addr^-(\sigma^m(\s))\neq\emptyset$ for all
   $m\geq 0$. Then 
   we can find a (possibly constant) sequence
   $\r^n\in \Addr^-(\s)$ with the
   property that each $\r^n$ is not an intermediate external address
   of length $\leq n+1$. By Corollary \ref{cor:boundeddistance}, 
   we may assume 
   that
   each $\r^n$ furthermore satisfies 
    \[ |r^n_j - s_j|\leq 1 \]
   for $j\in\{1,\dots,n\}$. It follows that every accumulation point of
   the sequence $\r^n$ is exponentially bounded. Since 
   $\Addr^-(\s)$ is compact, this concludes the proof. \end{proof}

 \begin{cor}[Rays accumulating on escaping points]
   \label{cor:raysaccumulating}
  Let $\Ek$ be an exponential map. Then one of the following holds:
  \begin{enumerate}
   \item There is some $\s\in\Sequ_0$ with
         $\Addr^-(\s)=\{\infty\}$; in particular, there are
         uncountably many addresses $\r\in\Sequ_0$ satisfying
         $\Addr^-(\r)=\{\r\}$; or
   \item every $\r\in\Sequ_0$ satisfies either
          $\Addr^-(\r)=\emptyset$ or $\Addr^-(\r)\cap\Sequ_0\neq\emptyset$. 
      \label{item:allexpbounded}
  \end{enumerate}
 \end{cor}
 \begin{proof} By the previous lemma, we see that (\ref{item:allexpbounded})
   holds unless there is some $\s\in\Sequ_0$ for which 
   $\Addr^-(\s)\neq\emptyset$ but $\Addr^-(\sigma(\s))=\emptyset$.
   This means that $\Addr^-(\s)=\{\infty\}$, as required. 
 \end{proof}

 \begin{proof}[Proof of Theorem \ref{thm:accumulation}]
   Using Lemma \ref{lem:accumulation}, 
   the Theorem is an immediate consequence of
   the previous Corollary and
   Theorem \ref{thm:raysatinfty}.
 \end{proof}

 We conclude this section by mentioning a few further questions suggested
  by this line of investigation.
 \begin{enumerate}
  \item If some dynamic ray accumulates at $\infty$, is it true that
    some ray must accumulate on the singular value?
  \item Can we replace $\Addr^-(\r)\cap\Sequ_0\neq\emptyset$ by
        $\r\in\Addr^-(\r)$ in Corollary \ref{cor:raysaccumulating}
        (\ref{item:allexpbounded})?
  \item Can we replace $\Addr^-(\r)\cap\Sequ_0\neq\emptyset$ by
        $\Addr^-(\r)=\{\r\}$ in Corollary \ref{cor:raysaccumulating}
        (\ref{item:allexpbounded})?
  \item If some dynamic ray accumulates on the singular value, can we
         ensure that there is some $\r$ for which $L_{\r}$ is an
         indecomposable continuum?
  \item If $\kappa\in J(\Ek)$, is there always some $\s\in\Sequ_0$ with
         $\Addr^-(\s)\neq\emptyset$?
 \end{enumerate}
 
 It seems reasonable to expect that the first two questions could be 
  answered using a further development of the methods in this section,
  while the remaining problems appear more difficult.

 \appendix

 \section{Remarks on higher-dimensional parameter spaces}
 
 For many entire functions with a bounded set of singular values, the
  escaping set consists of dynamic rays, just as in the exponential 
  family. In fact, this is now known to be true for
  all finite-order entire functions with a bounded set of
  singular values \cite{guenterthesis,fatoueremenko}. Let us shortly
  discuss how our main result generalizes to these cases.

 Schleicher \cite{dierkcosine} showed that, for postcritically 
  pre-periodic
  maps in the cosine family $z\mapsto a\exp(z)+b\exp(-z)$, every
  dynamic ray lands and every point of $\C$ is either on a ray
  or the landing point of a ray. Thus our main result is not
  true for the cosine family.

 The reason for this is that cosine maps do not have any asymptotic 
  values, and the presence of a dynamic ray landing at an asymptotic value 
  was the driving factor in our proof. Indeed, this is the only 
  obstruction to carrying our proof over to the general case. In
  particular, one would expect the following dichotomy: if 
  $f:\C\to\C$ is a postcritically finite entire function of finite order,
  then
  \begin{itemize}
   \item if $f$ has an asymptotic value, some dynamic ray of $f$ 
     accumulates on an entire dynamic ray, and conversely,
   \item if $f$ has only critical values, then every dynamic ray of $f$
     lands and every point in $J(f)$ is on a dynamic ray or the
     landing point of such a ray. 
  \end{itemize}

 We should note that the condition of ``accessible singular values''
  is also of interest in the study of these more general families of
  entire functions. In particular, under such a hypothesis it is
  possible to prove that all periodic dynamic rays land.
  (The only known
   proof that periodic rays of \emph{all} exponential maps land
   \cite{landing2new},
   regardless of whether the singular value is accessible,
   relies strongly on the one-dimensionality of
   exponential parameter space and does not generalize to the
   higher-dimensional case.)

\nocite{jackdynamicsthird}
\bibliographystyle{hamsplain}
\bibliography{/Latex/Biblio/biblio}

\end{document}